\documentclass[11pt]{article}
\usepackage{mathptmx}
\usepackage{amssymb}
\usepackage{latexsym}
\newtheorem{thm}{Theorem}[section]

\newtheorem{prop}[thm]{Proposition} \newtheorem{lemma}[thm]{Lemma}
\newtheorem{cor}[thm]{Corollary} \newtheorem{dfn}[thm]{Definition}
\newtheorem{conj}[thm]{Conjecture}
 \newtheorem{rmk}[thm]{Remark}
\newtheorem{ex}[thm]{Example} 

\newcommand {\pf}{\noindent{\bf Proof.}\ }

\newcommand{\reals}{{\mathbb R}}

\newcommand{\integers}{{\mathbb Z}}

\newcommand{\vol}{{\rm vol}}
\newcommand{\cok}{{\rm coker}}

\newcommand{\cald}{{\cal D}}
\newcommand{\cale}{{\cal E}}

\newcommand{\calh}{{\cal H}}

\newcommand{\calo}{{\cal O}}

\newcommand{\calv}{{\cal V}}
\newcommand{\calu}{{\cal U}}

\newcommand{\qed}{\begin{flushright} $\Box$\ \ \ \ \ \end{flushright}}

\newcommand{\frakh}{\mathfrak{h}}
\newcommand{\frakk}{\mathfrak{k}}

\newcommand{\frakt}{\mathfrak{t}}
\newcommand{\boe}{\mathbf{e}}
\newcommand{\bof}{\mathbf{f}}
\newcommand{\bog}{\mathbf{g}}

\newcommand{\arrows}{\,\lower1pt\hbox{$\longrightarrow$}\hskip-.24in\raise2pt
             \hbox{$\longrightarrow$}\,}
\newcommand{\defequal}{\stackrel{\mbox {\tiny {def}}}{=}}

\newcommand{\wedgetop}{\bigwedge^{\rm top}}
\begin{document}

\title{{\bf The volume of a differentiable stack}}
\author
{Alan
Weinstein
\thanks{Research partially supported by NSF Grant
DMS-0204100
\newline \mbox{~~~~}MSC2000 Subject Classification Number: 58H05
(Primary), 53D17 (Secondary).
\newline \mbox{~~~~}Keywords: Lie groupoid, Lie algebroid, modular
class, differentiable stack}\\
Department of Mathematics\\ University of California\\ Berkeley, CA
94720 USA\\ {\small(alanw@math.berkeley.edu)}}
\maketitle

\begin{abstract}
We extend the notion of the cardinality of a discrete groupoid
(equal to the Euler characteristic of the corresponding discrete
orbifold) to the setting of Lie groupoids.  Since this quantity is an invariant
under equivalence of groupoids, we call it the
volume of the associated stack
rather than of the groupoid itself.   Since there
is no natural measure in the smooth case
like the counting measure in the discrete case,
we need extra data to define the volume.  This data has the form of an
invariant section of a natural line bundle over the base of the
groupoid .   Invariant sections
of a square root of this line bundle constitute an ``intrinsic Hilbert
space'' of the stack. 

\end{abstract}

\section{Introduction}
\label{sec-intro}
As part of a larger program of groupoidification, 
the cardinality of a groupoid $G\arrows G_0$ with 
finitely many orbits and finite
isotropy groups was defined by Baez and Dolan
\cite{ba-do:finite} to
be the sum over orbits of the reciprocal of the orders of the isotropy
groups.    This sum
 is well defined because the isotropy groups of different
elements of an orbit are isomorphic.  It represents
the total mass of the
``natural measure'' of such a groupoid defined by Kim
\cite{ki:lefschetz} in connection with a Lefschetz formula for
equivariant cohomology.  This in turn was inspired by the same
expression in Behrend's Lefschetz formula \cite{be:lefschetz} for the
Frobenius automorphism on algebraic stacks.  The expression also
appears as the Euler characteristic of a 
0-dimensional orbifold.

 Baez and Dolan give
 many examples of and
reasons for their definition.  For instance, if the groupoid is the
transformation groupoid associated with the action of a group $H$ on a
space $X$, its cardinality is the quotient $\#(X)/\#(H)$ of the number
of elements 
$\#(X)$ in $X$ by the order $\#(H)$ of $H$.  In particular, if $X$ is
a single point, the 
cardinality is $1/\#(H)$.  

Since the cardinality is clearly an invariant under equivalence of
groupoids, we prefer to think of it as an invariant $\#(G_0//G)$ 
of the quotient stack\footnote{We refer to 
\cite{be:cohomology} and 
\cite{be-xu:differentiable}
for background material on differentiable stacks.} 
 $G_0//G$ (or $X//H$ in the case of a
transformation groupoid). 
The terminology is consistent with the usual one when the action is free,
so that $X//H$ is simply the (stack associated to the) set $X/H$, and 
$\#(X/H) = \#(X)/\#(H).$
On the other hand, when $X$ is a point, $X//H$ is the universal classifying
stack $BH$, so $\#(BH)$ becomes $1/\#(H)$.  As Baez and Dolan
themselves note, this result is
consistent with the idea that $BH={\rm pt}//H$ may be thought of as 
``one $\#(H)$'th of a point''.  (Do not try to read this aloud!)

The aim of the present work is to extend the notion of groupoid cardinality
from the discrete to the differentiable setting, i.e. to Lie groupoids
and their associated smooth stacks.  We call our extended notion
the {\bf volume of a
  stack}, thinking of the cardinality as a geometric rather than a
topological quantity.
  It is clear that we now need
additional data, namely measures which generalize the counting
measures on sets and groups which are implicit in the discrete situation.
It turns out that the appropriate data
are all contained in a $G$-invariant section of the bundle  
$Q_{A} \defequal \bigwedge^{\rm top}A\otimes\bigwedge^{\rm
  top}T^{*}G_0$ defined in \cite{ev-lu-we:transverse}, where
$A$ is the Lie algebroid of $G$.
  In fact, when $G=G_0$, such a section is just a volume
element $b$ on $G_0$, while when $G_0$ is a point, the section is the
value at the identity of a
bi-invariant multivector field of top degree, or ``inverse volume
element'', $a^{-1}$ (the inverse of the volume element $a$)
 on the group $G$.  In the
former case, our definition will give the integral\footnote{We will
  assume here and usually elsewhere that all our manifolds are
  oriented; the nonorientable case can be handled with the use of
  densities instead of forms.}
of $b$ over $G_0$, i.e. the volume of $G_0$ with the measure given by
$b$; in the latter, our definition will give the reciprocal of the
integral over $G$ of $a$.  (If
$a = 0$, or if $G$ is noncompact, the invariant is defined to be 0.)

We will see that, in the discrete case, our invariant reduces to that
of Baez and Dolan, 
with the additional flexibility that we can replace the counting
measures by arbitrary measures subject to an invariance condition.
Furthermore, if $G$ is the groupoid associated with a smooth action of a
Lie group $H$ on a manifold $X$, a special kind of invariant section of
$Q_A$ is the product of an $H$-invariant volume element $b$ on $X$ and a
bi-invariant, nowhere-vanishing volume element $a$ on $H$, and for such a section 
our volume for $X//H$ will be the quotient $\int_X b / \int_H a.$
This quotient makes clear sense when $X$
and $H$ are compact; when this is not the case, but the action of $H$
on $X$ is still proper and cocompact (i.e. having compact orbit
space\footnote{Perhaps another appropriate name for this condition
  would be ``coproper''.}),
our invariant continues to be well-defined.  Note that $H$ is not
required to be unimodular; when it is not, the non-invariance of a Haar measure
under conjugation must be countered exactly by the non-invariance of a
measure on $X$.

More generally, we will see that any invariant section of $Q_A$
defines a measure on the coarse moduli (or ordinary quotient, or
orbit) space
$G_0/G$ whose integral over any relatively compact open subset $U$ is the
volume of the inverse image of $U$ under the natural projection from
$G_0//G$ to $G_0/G$.  Another way of expressing this is to say that 
a section of
$Q_A$ {\em is} 
a ``smooth measure'' on the stack, and the ordinary measure on the
orbit space is its push-forward.   We can also consider the line
bundle associated to the frame bundle of 
$Q_A$ via the homomorphism $a\to |a|^{1/2}$ from
the multiplicative reals to themselves.  As in \cite{we:maslovgerbe},
we denote this line bundle by $|Q_A|^{1/2}$.  $G$ still acts on this
bundle, and its space of 
compactly supported invariant sections carries a natural inner
product.  The Hilbert space completion of this space is an ``intrinsic 
$L^2$ space'' for the stack $G_0//G$.  In a similar way, one may
construct $L^p$ spaces in duality, and spaces of distributions.

A natural next step from here should be to attach a 
vector space (sections, or distributional sections of $Q_A$ or one of
its powers) to each groupoid $G$, and to 
attach a linear map to each morphism in some
nice category.  Here, we think of the volume of a stack as a linear
map $\reals \to \reals$ associated to the diagram of groupoids ${\rm
  pt}\leftarrow G \rightarrow {\rm pt}.$  Such a diagram is a special
case of a {\bf span} $G \leftarrow S \rightarrow H$, which, according
to the degroupoidification
program described in \cite{ba-ho-wa:groupoidification}, 
should be considered as a
kind of morphism from $G$ to $H$.  

It may be that the appropriate morphisms to consider when extending
from the finite case should be those which are between stacks of the
same dimension, or perhaps submersions of stacks, for which one could
try to integrate over the fibres.  Here, an extension to groupoids of
the 
relative modular class for Lie algebroid morphisms in  
\cite{ko-la-we:modular} should come into play.

From the point of view of microlocal analysis, one might even want to
consider more general geometric morphisms between stacks, encoded by lagrangian
submanifolds and symbols, which would induce pseudodifferential or
even Fourier integral operators between the corresponding vector
spaces.  A
prerequisite for doing this would be a correct definition of the
cotangent bundle of a stack as a ``symplectic stack''.  A symplectic
stack cannot be simply a stack defined by a symplectic
groupoid, since such a stack has no symplectic structure of its own.  
Rather, the cotangent stack should be presented by a ``groupoid in the
symplectic category'', in which the groupoid structure operations are canonical
relations which may even be multiply-defined.   In fact, one 
gets precisely such an object when applying to the structure
operations in a Lie groupoid the functor which assigns to each
mapping (or relation) between manifolds the conormal bundle to its
graph.  
 
\bigskip
\noindent
{\bf Acknowledgments} I would like to think John Baez, Kai Behrend, Rui Loja
Fernandes, Minhyong Kim, Yvette Kosmann-Schwarzbach, Eckhard
Meinrenken, Martin Olsson for helpful discussion and comments.   
I would also like to thank the group Analyse Alg\'ebrique at the
Institut Math\'ematique de Jussieu for their hospitality.

\section{First definition and the finite case}
In this section, we give a provisional definition of the volume
of the stack $G_0//G$ (sometimes denoted $BG$) presented by a groupoid
$G\arrows G_0$ (sometimes denoted simply as $G$)
in terms of an invariant 
section $\lambda$ of $Q_A$, where $A$ is the Lie algebroid of $G$.
We will see that this definition has two difficulties.  First of all, it 
depends on the
decomposition of $\lambda$ as the quotient of a section $b$ of 
$\bigwedge^{\rm top}T^*G_0$ by a section $a$ of $\bigwedge^{\rm
  top}A^*$, and it is not clear that it is independent of the
decomposition.   Second, it produces a sensible result only when the
source and target maps from $G$ to $G_0$ are proper, although only
properness of the groupoid, i.e. properness of the map
$(l,r):G \to G_0 \times G_0$, should be needed.   (Here and later, 
the letters $l$ and
$r$ will denote the source and target maps from $G$ to $G_0$.  The
reader may choose which is to be the source and which is the target,
but we will insist that the product $gh$ is defined when
$r(g)=l(h)$.)

In the next section, we will see how to remove the difficulties
with our definition.  But before even giving the definition, we will 
look more closely at the finite case.  

\begin{dfn}[Baez-Dolan\cite{ba-do:finite}]\label{def:finite}
If $G=G\arrows G_0$ is a groupoid with finite isotropy groups and
finitely many orbits, we define the {\bf cardinality of the stack}
$G_0//G$ to be the sum 
$$\#(G_0//G) \defequal \sum_{\calo} \#(G_{\calo})^{-1},$$ where
$\calo$ ranges over the orbit space $G_0/G$, and $\#(G_\calo)$ denotes
the cardinality of the isotropy of any element of $\calo$.  (The
cardinality of the empty groupoid is zero.)
\end{dfn}

It is clear from the definition that equivalent groupoids give rise to
the same cardinality, so that this quantity is really an invariant of the
stack, independent of its presentation by a groupoid.  

As Baez and Dolan \cite{ba-do:finite} note, 
the corresponding sum may be convergent for
some interesting groupoids with infinitely many orbits.  For example,
the cardinality of the stack presented by the groupoid of all finite
sets and their bijections (or an equivalent groupoid whose objects
form a set rather than a class) is $e = 2.718\ldots$.  

It is not so simple to transfer Definition \ref{def:finite} to the
smooth case, since the orbit space is generally not a smooth manifold.
Instead, we will reformulate the definition on the basis of the 
following simple fact, which is a
special case of Theorem \ref{thm:finite} below.  

\begin{prop}\label{prop:finite}
For any finite groupoid $G\arrows G_0$,  
$$\#(G_0//G) = \sum_{y\in G_0}\#(r^{-1}(y))^{-1}.$$
\end{prop}

Proposition \ref{prop:finite} suggests defining the volume of
$G_0//G$, when $G$ is a Lie groupoid, by integrating over $G_0$ the
reciprocals of the volumes of the fibres of the surjective submersion
$r:G\to G_0$.  To do this, we need a measure on $G_0$ and measures
along the fibres of $r$.  Using this data, we make the following
provisional definition, assuming as usual that $G$ and $G_0$ are oriented.  

\begin{dfn}\label{firstdefinition}
Let $G\arrows G_0$ be a compact Lie groupoid (i.e.~$G$ is a compact
manifold, so that $G_0$ and the fibres of $r$ are also compact) with
Lie algebroid $A$.
Let $a$ be a nowhere vanishing section of $\wedgetop A^*$ and $b$ a
section of $\wedgetop T^*G_0$.  The {\bf volume of the stack} $G_0//G$ 
with data $(a,b)$ is defined as
$$\vol _{(a,b)}(G_0//G) \defequal \int_{y\in G_0}\left(\int_{r^{-1}(y)}a_r\right)^{-1} b . $$
\end{dfn}

In this definition, $a_r$ is the right-invariant form, defined along the tangent
bundle to the $r$-fibres, whose values along the unit section are
given by $a$.  (The Lie algebroid $A$ is identified with the tangent
bundle along the units to the $r$-fibres.)

To show that this definition depends only on the stack and not on the
presenting groupoid, we cannot even begin without having
 a way of moving the data $(a,b)$
from one groupoid to any equivalent one.  But this is not possible;
what is transferable between equivalent groupoids is only the product
$a^{-1}b$, which is a section of the tensor product line bundle 
$Q_A \defequal \wedgetop A \otimes \wedgetop T^*G_0.$  
(See \cite{ko-la-we:modular} and Section \ref{sec-morita} below.)
We thus need to show that $\vol _{(a,b)}G_0//G$
depends only the product $a^{-1} b$ and not on the individual factors;
i.e. that the volume is unchanged when $a$ and $b$ are multiplied by
the same nonvanishing function $\theta$ on $G_0$.  

To see how to proceed, we return to the finite case, where the bundles
$A$ and $TG_0$
 have 0-dimensional fibres, so sections of their top exterior powers
(and of the top exterior powers of their duals) 
are simply scalar functions.  With data $(a,b)$, Definition
\ref{firstdefinition} in the case of a finite groupoid becomes
$$\vol_{(a,b)}(G_0//G) \defequal \sum_{y\in G_0}\left(
\sum_{g\in r^{-1}(y)} a(l(g))\right) ^{-1} b(y),$$
which  becomes the cardinality formula in Proposition 
\ref{prop:finite} when $a$ and $b$ are unity. 

For the outer sum, we may sum over each orbit and then sum over
the orbit space, i.e.
$$\vol_{(a,b)}(G_0//G) = \sum_{\calo \in G_0/G} \, \sum_{y\in \calo}\left(
\sum_{g\in r^{-1}(y)} a(l(g))\right) ^{-1} b(y).$$
For each orbit $\calo$, we have
$$S_\calo \defequal \sum_{y\in \calo}\left(
\sum_{g\in r^{-1}(y)} a(l(g))\right) ^{-1} b(y)=
\sum_{y\in \calo}\left(\sum_{x\in \calo}\sum_{g\in l^{-1}(x)\cap r^{-1}(y)}
  a(x)\right)^{-1} b(y).  
$$
Now the number of elements in $l^{-1}(x)\cap r^{-1}(y)$ depends only
on the orbit $\calo$ and is equal to the cardinality $\#(G_\calo)$ of
the typical isotropy group.  Hence, we have
$$S_\calo
=  \#(G_\calo)^{-1} \left( \sum_{x\in \calo} a(x) \right)^{-1} \sum_{y\in \calo} b(y).  
$$

In general, $S_\calo$ will depend on all the values of $b/a$ on $\calo$,
but if $\lambda=b/a$ is constant on orbits (i.e. a $G$-invariant section of
the trivial bundle $Q_A$), it only depends on the constant value
$\lambda(\calo)$ of that section on the orbit, and we obtain the final
formula:
$$\vol_\lambda (G_0//G) = \sum_{\calo\in G_0/G} \#(G_\calo)^{-1} \lambda(\calo).$$
Since the right hand side is clearly invariant under equivalence of
groupoids, we have the following result.

\begin{thm}\label{thm:finite}
Let $G\arrows G_0$ be a finite
groupoid.  Let $a$ and $b$ be functions on $G_0$ such that 
$a$ is nowhere vanishing and the quotient $\lambda = b/a$ is
$G$-invariant, so that $\lambda$ may be considered as a function on $G_0/G$,
or a $G$-invariant section of $Q_A$, where $A$ is the
(zero-dimensional) Lie algebroid of $G$.

Then the quantity $$\sum_{y\in G_0}\left(
\sum_{g\in r^{-1}(y)} a(l(g))\right) ^{-1} b(y)$$
is equal to 
$$\sum_{\calo\in G_0/G} \#(G_\calo)^{-1} \lambda(\calo).$$
In particular, it depends on $a$ and $b$ only via their quotient
$\lambda$.  

Furthermore, given an equivalence between $G$ and another finite groupoid
$G'$, with Lie algebroid $A'$, there is a bijective correspondence 
between $G$-invariant
sections of $Q_A$ and $Q_{A'}$, and if $\lambda ' = b'/a'$ is the section corresponding
to $\lambda = b/a$, then
$$\sum_{y\in G_0}\left(
\sum_{g\in r^{-1}(y)} a(l(g))\right) ^{-1} b(y) = 
\sum_{y'\in G_0 '}\left(
\sum_{g'\in r'^{-1}(y')} a'(l'(g'))\right) ^{-1} b'(y').$$
\end{thm}

We may therefore make the following definition.

\begin{dfn}
Let $G\arrows G_0$ be a groupoid with finite isotropy groups, $G_0//G$
the corresponding stack, and $\pi:G_0//G \to G_0/G$ the natural projection.
Then any function $\lambda$ on $G_0/G$ (i.e. $G$-invariant section of
$Q_A$, where $A$ is the Lie algebroid of $G$) defines a measure
$\mu_{\lambda}$ on
$G_0//G$ for which the ``measurable sets'' are the preimages under
$\pi$ of finite subsets of $G_0/G$, and the measure of
such a subset $\calu$ is defined by
$$\mu_{\lambda}(\pi^{-1}(\calu))\defequal \sum_{\calo\in \calu}
\#(G_\calo)^{-1} \lambda(\calo).$$  
\end{dfn} 

We may therefore consider $\mu_\lambda$ 
as a volume element, or measure, on $G_0//G$;  
its push-forward under the natural projection
$\pi:G_0//G\to G_0/G$ is the measure which assigns to each point $\calo$
the measure $\#( G_\calo) ^{-1} \lambda (\calo)$.  

\section{The differentiable case}
Definition \ref{firstdefinition} has the virtue that it clearly leads
to a finite result for any compact groupoid, but its invariance
properties are hard to verify directly.  We would like to imitate
the orbit decomposition method of the previous section, but this works nicely
only when the groupoid is strongly 
regular in the sense that the decomposition into
orbits is a fibration.  Our strategy will be to apply the orbit
decomposition on the strongly regular part of the groupoid, whose complement
turns out to be negligible as far as integration is concerned.  The
latter fact follows immediately from the slice theorem in \cite{we:proper}
and Zung's linearization theorem \cite{zu:proper}  for proper groupoids, 
since the orbit structure of
a proper groupoid is locally like that of the action of a compact
group, for which there is a principal orbit type.

Assume now, then, that $G\arrows G_0$ is a strongly regular, compact
(hence proper) groupoid, and let $f:G_0\to G_0/G$ be the natural
projection.  Then the integral over $G_0$ in Definition
\ref{firstdefinition} of the function 
$$\left(\int_{r^{-1}(y)}a_r\right)^{-1}$$ times the volume element $b$
can be written as iterated integral--first over
the fibres of $f$, i.e. over the individual orbits, and then over the
orbit space.  To do this, we will need to decompose $b$ as the product of a 
volume element along the $G$-orbits and one on the orbit
space .  

The inner integral in Definition \ref{firstdefinition} can also be
written as an iterated integral.  In fact, the restriction to
$r^{-1}(y)$ of the map $l$ is a principal fibration over the orbit
$Gy$ with structure group the isotropy 
group $G_y$ acting from the right by groupoid multiplication.  Again,
to get an iterated integral, we must decompose the integrand $a_r$ as
a product of a volume element along the $G_y$-orbits and one along the
$G$-orbit $Gy$.  

To obtain the decompositions above we must make a choice.  To see what
to do, we recall the exact sequence of vector bundles over $G_0$,
$$0\to \ker\rho \to A \to TG_0 \to \cok\rho\to 0,$$ where $\rho:A\to
TG_0$ is the anchor map of $G$, which may be identified with the
restriction of $Tl:TG\to TG_0$ to the kernel of $Tr$
along the unit section.  The kernel of
$\rho$ is thus the bundle of Lie algebras of the isotropy groups,
while its cokernel is the conormal bundle to the 
foliation by $G$-orbits.   The standard ``alternating product'' rule
for top exterior powers in an exact sequence yields a natural
isomorphism 
\begin{equation}\label{eq:naturalisomo}
Q_{A}=\bigwedge^{\rm top}A
\otimes\bigwedge^{\rm  top}T^{*}G_0 \approx \bigwedge^{\rm
  top}\ker\rho\otimes\bigwedge^{\rm top} (\cok\rho)^*. 
\end{equation}
More explicitly, the isomorphism in (\ref{eq:naturalisomo})
   comes from the natural isomorphisms
$A/\ker \rho \approx \rho(A)$ and $\cok\rho \approx TG_0/\rho(A).$

Now, given  data $a$ and $b$ as above for which $\lambda = a^{-1} b$ is a 
$G$-invariant section of $Q_A$, there is a corresponding $G$-invariant
section, which we also denote by $\lambda$, of $\bigwedge^{\rm
  top}\ker\rho\otimes\bigwedge^{\rm top} (\cok\rho)^*$.  Since $G$
is a proper groupoid acting on the line bundle 
$\bigwedge^{\rm  top}\ker\rho$ (which we assume, as usual, to be
orientable), we can find
a
non-vanishing invariant section\footnote{To construct such a section,
start with any non-vanishing section, then average over $G$, using a
cutoff function \cite{tu:hyperboliques} as in the proof of the
vanishing theorem for groupoid cohomology in \cite{cr:differentiable}.}
 $\alpha$ of this bundle. We may then 
write $\lambda$ as $\alpha^{-1} \beta$, and the invariance of
$\lambda$ implies that $\beta$ is an invariant section of 
of $\bigwedge^{\rm top} (\cok\rho)^*$.  

\begin{rmk}
\upshape
Note that we could not have
imposed the invariance requirement on $a$ and $b$, since the line bundles of
which they are sections do  not have natural actions of $G$, only 
actions up to homotopy.   On the other hand, it is only
in the regular case that we can speak of smooth sections of the kernel
and cokernel of $\rho$.   
\end{rmk}

The section $\beta$ gives a bi-invariant volume element on each
isotropy group of 
$G$ which is ``the same'' on all the isotropy groups over a given
orbit.  We thus have a well-defined volume function $\vol_\beta
(G_\calo)$ on the orbit space $G_0/G$.  On the other hand,
$\alpha$ is a $G$-invariant volume element on the orbit space.  We
therefore have the well defined expression
$$\int_{\calo\in G_0/G} \vol_\beta(G_\calo)^{-1} \alpha,$$ 
which reduces to 
$$\sum_{\calo\in G_0/G} \#(G_\calo)^{-1} \lambda(\calo)$$
in the finite case.  
 Furthermore, it is
  clear that this expression does not depend on the choice of $\beta$,
  since multiplying it by a function $\theta$, which must be invariant, requires multiplication
  of $\alpha$ by the same function, and the effects of the two
  multiplications cancel one another.

We will now prove the following extension of Theorem \ref{thm:finite},
which validates Definition \ref{firstdefinition}

\begin{thm}\label{thm:smooth}
Let $G\arrows G_0$ be a compact Lie 
groupoid with Lie algebroid $A$ for which the map $r:G\to G_0$ is 
proper (hence a locally trivial fibration).  Let $a$ and $b$ be sections of 
$\bigwedge^{\rm top}A^*$ and $\bigwedge^{\rm top}T^{*}G_0$ respectively such that 
$a$ is nowhere vanishing, and the quotient $\lambda = a^{-1}b$ is a
$G$-invariant section of $Q_{A} \defequal \bigwedge^{\rm
  top}A\otimes\bigwedge^{\rm top}T^{*}G_0$.

Then the quantity $$\int_{y\in G_0}\left(
\int_{r^{-1}(y)} a_r\right) ^{-1} b$$
is equal to
$$\int_{\calo\in G_0/G} \vol_\beta(G_\calo)^{-1}\alpha,  $$
where $\alpha$ and $\beta$ are any $G$-invariant sections of 
$\wedgetop (\ker\rho)^*$ and $\wedgetop (\cok \rho)^*$ respectively
such that $\alpha$ is nowhere vanishing and $\alpha^{-1} \beta$
corresponds to $\lambda$ under the natural isomorphism
(\ref{eq:naturalisomo}).  In particular, it depends on $a$ and $b$
only through their quotient $\lambda$.

Furthermore, given an equivalence between $G$ and another compact groupoid
$G'$, with Lie algebroid $A'$, there is a 
bijective correspondence between $G$-invariant
sections of $Q_A$ and $Q_{A'}$, and if $\lambda ' = b'/a'$ is the section corresponding
to $\lambda = b/a$, then
$$\int_{y\in G_0}\left(
\int_{r^{-1}(y)} a_r \right) ^{-1} b = 
\int_{y'\in G_0 '}\left(
\int_{r'^{-1}(y')} a'_{l'}\right) ^{-1} b'.$$
\end{thm}

\pf
Since a compact Lie groupoid is proper, and a proper groupoid is
locally equivalent to action groupoids for actions of compact groups,
the natural projection $\pi:G_0\to G_0/G$ is a fibration with compact fibres
when restricted to an invariant open subset $\calv \subset G_0$ whose
complement has positive codimension and therefore does not contribute
to the integral in the theorem.  We may therefore take the integral
over this strongly regular set and will use the Fubini theorem, imitating the
summation proof of Theorem \ref{thm:finite}.   We denote by $\omega$
the unique section of $\wedgetop \rho(A)^*$ for which $b=\omega\beta$
(from which it follows that $a=\alpha\omega$).  
\begin{eqnarray*}
\int_{y\in G_0}\left(
\int_{g\in r^{-1}(y)} a_r\right) ^{-1} b 
& =  & \int_{y\in \calv}\left(\int_{r^{-1}(y)} a_r\right) ^{-1} b 
\\ &=& \int_{y\in \calv}\left(\int_{r^{-1}(y)} a_r\right) ^{-1} \omega\beta
\\ &=& \int_{\calo\in \pi(\calv)}
\left[ \int_{y\in\calo}\left(\int_{r^{-1}(y)} (\alpha\omega)_l\right) ^{-1} 
 \omega\right]\beta.
\end{eqnarray*}
But, using the fibration $l:r^{-1}(y)\to\calo$, we have
\begin{eqnarray*}
\int_{r^{-1}(y)}(\alpha\omega)_l 
& = & \int_{x\in\calo} \left( \int_{l^{-1}(x)\cap    r^{-1}(y)}\alpha\right)\omega
\\ &=& \int_{x\in\calo} \vol_\alpha(G_\calo) \omega
\\ &=& \vol_\alpha (G_\calo) \int_\calo \omega.
\end{eqnarray*}
Combining the last two calculations, we obtain:
\begin{eqnarray*}
\int_{y\in G_0}\left(
\int_{g\in r^{-1}(y)} a_r\right) ^{-1} b 
& = & \int_{ \calo \in
  \pi(\calv)}\left[\int_{y\in\calo}\left(\vol_\alpha(G_\calo)\int_\calo
    \omega\right)^{-1} \omega\right]\beta
\\& = & \int_{ \calo \in
  \pi(\calv)}\vol_\alpha(G_\calo)^{-1}\left(\int_\calo
    \omega\right)^{-1} \left(\int_\calo  \omega\right)\beta
\\& = & \int_{\calo\in \pi(\calv)}\vol_\alpha(G_\calo)^{-1}\beta
\\& = & \int_{G_0/G} \vol_\alpha(G_\calo)^{-1}\beta.
\end{eqnarray*}
(For the last equality, we simply take integration over the strongly regular
part as a definition of integration over the singular space $G_0/G$.
Convergence of the integral over the noncompact complement of the
singular points is guaranteed by its equality with the expression
involving nonsingular integrals over $G$ and the $r$-fibres.)
\qed

We may therefore make the following definition.
\begin{dfn}\label{seconddefinition}
Let $G\arrows G_0$ be a proper Lie groupoid, $G_0//G$
the corresponding stack, $\pi:G_0//G \to G_0/G$ the natural
projection, and $\calv \subseteq G_0$ the strongly regular set.
Then any $G$-invariant section $\lambda$ of
$Q_A$, where $A$ is the Lie algebroid of $G$, defines a (signed) measure
$\mu_{\lambda}$ on
$G_0//G$ for which the ``measurable sets'' are the preimages under
$\pi$ of relatively compact open subsets of $G_0/G$, and the measure of
such sets is defined as
$$\mu_{\lambda}(\pi^{-1}(\calu))\defequal \int_{\calo\in \calu\cap \pi(\calv)}
\vol_\beta(G_\calo)^{-1} \alpha,$$ where $\alpha$ and $\beta$ are as
in Theorem \ref{thm:smooth}.  
\end{dfn}

We may therefore consider $\lambda$ as a 
volume element, or measure
on the stack $G_0//G$.
Its push-forward under the natural projection
$\pi:G_0//G\to G_0/G$ is the measure which assigns to each relatively
compact
open subset $\calu$ of $G_0/G$ the integral of the 
smooth form $\vol_\beta(G_\calo)^{-1} \alpha$ over the intersection of
$\calu$ with the part of $G_0/G$ over which $G$ is strongly regular.  

\section{Morita invariance}
\label{sec-morita}
Although the formula in Definition \ref{seconddefinition} shows that
the volume of a stack clearly depends only on Morita invariant data
for a groupoid, it still needs to be shown that
any Morita equivalence between groupoids $G'$ and
$G''$ with Lie algebroids $A'$ and $A''$ induces a natural isomorphism
between the spaces of smooth invariant sections of the line bundles
$Q_{A'}$ and $Q_{A''}$.  These spaces may therefore be considered as
representatives of ``the space of volume elements of a stack''. 

The argument for invariance begins with the special case of the
inclusion of an open subset of $G_0$ which is {\bf full}, i.e which 
 intersects every
$G$-orbit. 

\begin{lemma}
\label{lem-inclusion}
Let $G\arrows G_0$ be a groupoid, $G'_0\subseteq G_0$ a full open
subset, $G'\arrows G'_0$ the restriction of $G$ 
to $G'_0$.    
Then the restriction operation defines an isomorphism from
the smooth $G$-invariant sections of $Q_A$ to the smooth 
$G'$-invariant sections of $Q_{A'}$.  
\end{lemma}
\pf 
Since $G'_0$ intersects each orbit, each invariant section
$\lambda'$ of $Q_{A'}$ extends uniquely to an invariant section of
$Q_A$.  It is clear that restriction preserves smoothness; we just
have to show the same for extension, i.e. that an invariant section of
$Q_A$ is smooth if and only if its restriction to some full open subset
is smooth.  

It is easy to prove the last statement when $G_0$ has constant
dimension, using smooth bisections, which exist through each point of
$G$.  But we will need to apply the result to the case where $G_0$ has
components of different dimensions, and so we take another approach.  

Note that the action of $G$ on $Q_A$ can be encoded as an isomorphism
$j:l^! Q_A \to r^! Q_A$ of line bundles over $G$.  
Following Mackenzie \cite{ma:general}
we use $!$ rather than $*$ to denote
pull-backs.)  A section $\lambda$ of $Q_A$ 
is then $G$-invariant if and only if the pulled back sections satisfy
$j l^! \lambda = r^! \lambda$.  Moreover, since $l$ and $r$ are submersions, 
$\lambda$ is smooth in a
neighborhood of $x \in G_0$ if and only
if the pull-back $l^! \lambda$ (or 
$r^! \lambda$) is smooth in a neighborhood of some point of $l^{-1}(x)$
(or $r^{-1}(x)$).  It follows easily that the set of points
$x\in G_0$ near which an invariant section $\lambda$ is smooth
 is $G$-invariant, and the lemma follows.
\qed

We thus have:

\begin{cor}
\label{cor-twoinclusions}
Let $G\arrows G_0$ be a groupoid, $G'_0$ and $G''_0$ full open subsets.
Then the restrictions from $G_0$ induce a
natural isomorphism between smooth invariant 
sections of $Q_{A'}$ and those of $Q_{A''}$.
\end{cor}

To establish invariance of $Q_A$ under general Morita equivalences, 
we will use the following
standard\footnote{It is hard to pinpoint the first occurrence of this
  notion, but it is implicit in work of Kumjian and Renault giving
  geometric constructions of linking algebras in the algebraic theory of 
Morita equivalence.   The first explicit appearance of the term seems
to be in \cite{mu:bundles}.}
    notion of  ``linking groupoids''.  

\begin{prop}
\label{prop-linking}
Let $G'$ and $G''$ be groupoids and $B$ an invertible
$(G',G'')$-bibundle.  Then there is a unique (up to natural
isomorphism) groupoid $G$ such that:
\begin{enumerate}
\item  $G_0$ is the disjoint union of
$G'_0 $ and $G''_0$.
\item  $G'$ and $G''$ are the restrictions of $G$ to
$G'_0$ and $G''_0$ respectively.
\item  The set $l^{-1}(G'_0) \cap r^{-1}(G''_0)$ is identified with
  $B$ in such a way that the projections of $B$ to  $G'_0$ and $G''_0$ and the
left and right actions on the bibundle $B$ 
with left and right multiplication in
  the groupoid $G$.  
\end{enumerate}
\end{prop}
\pf
We assume without loss of generality that $B$ is disjoint from $G'$
and $G''$.  Let $G$ be the disjoint union $G'\cup B \cup
\overline{B}\cup G''$, where $\overline{B}$ is a copy of $B$.  The
inversion operation on $G$ is defined to be the union of those on $G'$
and $G''$ and the correspondence between $B$ and $\overline{B}$, along
with its inverse.  The groupoid structure on $G$ is then completely
determined by the conditions enumerated in the statement of the Proposition.
\qed

Combining Corollary \ref{cor-twoinclusions} and Proposition
\ref{prop-linking}, we obtain the main result of this section.  The
proof of the last part of the following statement is left to the reader.

\begin{thm}
\label{thm-morita}
Let Let $G'$ and $G''$ be groupoids and $B$ an invertible
$(G',G'')$-bibundle.   Then $B$ induces an isomorphism between the
$G$-invariant sections of $Q_{A'}$ and the $G'$-invariant sections of
$Q_{A''}.$   Isomorphic bibundles induce the same isomorphism, and the
composition of bibundles induces the composed isomorphism.
\end{thm}

\begin{rmk}
\upshape
A more formal way of stating the results above is that
there is a 2-functor from the 2-category of groupoids,
invertible bibundles, and isomorphisms of invertible bibundles to the
``discrete'' 2-category of vector spaces, isomorphisms of vector
spaces, and trivial isomorphisms of isomorphisms.  This is very
reminiscent of the decategorification construction of
\cite{ba-do:finite}, except that it is not clear how it extends
to more general morphisms.   It may be that the more appropriate 
construction is that which produces (sometimes partially defined)
morphisms of vector spaces of distributions from
suitable canonical relations between cotangent bundles, as suggested 
in the introduction.
\end{rmk}

\section{Transformation groupoids}

Let $h\mapsto h_X$ be an action
of a Lie group $H$ on a manifold $X$, and let $v\mapsto v_X$ be the
corresponding Lie algebra action, which takes $v\in \frakh$
to the vector field on $X$ which generates 
the 1-parameter group $(\exp (-tv))_H$ of diffeomorphisms.
(We need the minus
sign to get a representation rather than an antirepresentation.)
Let $G=H\times X \arrows X =
G_0$ be the corresponding transformation Lie groupoid.  Its Lie
algebroid is the trivial bundle $A=\frakh \times X$, and the anchor
is $(v,x)\mapsto v_X(x).$  Any section of $Q_A$ may be factored as
$a^{-1} b$, where $a$ is a constant element of $\wedgetop \frakh ^*$ and
$b$ is a form of top degree on $X$.  If $H$ is unimodular, we can take
$a$ to be adjoint-invariant, in which case $G$-invariance of $a^{-1}b$
is equivalent to $H$-invariance of $b$.  If $H$ is not unimodular,
then, for $a^{-1}b$ to be invariant, we must have, for each $h\in H$, 
$h_X^* b = \mu(h) b$, where $\mu:H\to \reals^{\times}$ is the modular
function of $H$ (i.e. the determinant of the adjoint representation).

Whether or not the group $H$ is unimodular, the transformation
groupoid is unimodular as long as the action (and hence the groupoid)
is proper, which we will assume from now on.  As a result, we can
choose $a$ and $b$ as required above.  Since the isotropy groups are
compact, the modular function is identically $1$ on them, and so the
form $b$ will be $H$-invariant.  

If $H$ is compact, we may apply Definition 
\ref{firstdefinition} to any relatively compact region in 
$X/H$ to conclude that the induced measure there, pushed forward from
$X//H$, is $\left(\int_H a_r\right)^{-1} b$.  In particular,
if the integral of $b$ over $X$ is finite, we have
$$\vol_{(a,b)}(X//H) = \vol_b(X)/\vol_a(H),
$$
exactly as in the finite case.

If $H$ is not compact, we must use Definition
\ref{seconddefinition}.  As in the discussion leading up to that
definition, we choose, over the regular part of $G_0$, an invariant
section $\alpha$ of $\wedgetop \ker \rho ^*$, i.e. an invariant family of
bi-invariant measures 
on the isotropy groups.  These, together with $a$, induce
a section of $\wedgetop(A/\ker \rho)^*.$  This section is transferred by
$\rho$ itself to a section of $\wedgetop \rho(A)^*$,
which with $b$ then induces an
invariant section $\beta$ of $\cok \rho$, i.e. a top-degree form on
the orbit space.  We may then use the formula in Definition
\ref{seconddefinition} to compute the measure on $X//H$ arising from
the invariant section $a^{-1}b$ of $Q_A$.  

\begin{ex}\label{ex:plane}
\upshape
Let $X$ be the euclidean plane with its usual euclidean measure
$b=dx\wedge dy = r
dr\wedge d\theta$, and
$H$ the group $SO(2)$ acting in the usual way by rotations around the
origin.  Let $a$ be the usual angle measure on $SO(2)$.  Then the
isotropy groups of regular points are trivial; if we take counting
measure on them, the induced measure
on the regular orbit space $\reals^+$ with coordinate $r$ is the
quotient $(d\theta)^{-1} \otimes r dr \wedge d\theta $, which may be
identified with $r dr$.  The measure on the stack $\reals^2//SO(2)$ is then
given by $r dr$.  

Suppose, now, that we replace $SO(2)$ by the full orthogonal group
$O(2)$.  Since the adjoint representation of $O(2)$ is no longer
orientation preserving, we must work with densities rather than forms,
but we may take the infinitesimal data $a^{-1} b$ to be essentially
the same as before.  The main difference here is that the isotropy
groups are now $\integers _2$, so that the measure on the stack
$\reals^2//O(2)$ becomes $\frac{1}{2} r dr$.  
  \end{ex}

\subsection{Stacks of conjugacy classes and adjoint orbits}

Let $H$ be any Lie group acting on either $H$ itself by conjugation or
on the Lie algebra $\frakh$ by the adjoint representation.  In either
case, the space $X$ on which $H$ is acting carries a natural
isomorphism from $\wedgetop T^*X$ to the trivial bundle with fibre
$\wedgetop \frakh^ *$.  As a result, the bundle $Q_A$ has a natural invariant
trivialization, so the groupoid is unimodular, and
we may use the constant section $1$ (i.e. choosing
$a$ and $b$ to be ``the same'') to compute ``canonical'' measures on $H//H$ and
$\frakh//H$ when the action is proper.  In fact the action is proper
just when $H$ is compact.   In this case we have $\vol_1 (H//H) = 1,$
while the volume of $\frakh//H$ is infinite.   In either case, it is
still interesting to compute the induced measure.  We concentrate
on the case of the action on $\frakh$, since that on $H$ is related to
it by the exponential map, which is equivariant.  (The regular part of
$H$ is contained in the set of regular values of the exponential map.)

To compute the induced measure, we begin by 
choosing invariant measures on the isotropy groups of the strongly
regular elements, i.e. on the maximal tori.   For convenience, we
choose the measures for which the total volume is $1$.  

Following \cite{do-wi:harmonic},  we 
choose a basis of the Lie algebra $\frakh$ of the form
$$(\boe_1,\ldots,\boe_r,\bof_1,\ldots,\bof_k,\bog_1,\ldots,\bog_k),$$
where the first $r$ entries are a basis of a Cartan subalgebra
$\frakt$ which are also a basis for the lattice
$\exp^{-1}(e) \cap \frakt$, 
and, for each $j$, $\bof_j$ and $\bog_j$ span a
plane on which the adjoint action of $\frakt$ is given by
$[\boe,\bof_j]=\sigma_j(\boe)\bog_j$ and  $[\boe,\bog_j]=-\sigma_j(\boe)\bof_j$ 
for linear forms $\sigma_1,\ldots,\sigma_k$ in $\frakt ^*$.  The dual
basis will be denoted
$(\boe_1^*,\ldots,\boe_r^*,\bof_1^*,\ldots,\bof_k^*,\bog_1^*,\ldots,\bog_k^*),$ 
As measure\footnote{We will abuse language by
  using the term ``measure'' to refer to top degree forms, and sometimes
  even to particular values of these forms.}   
 on $\frakt$ we choose $\boe_1^*\wedge\ldots\wedge\boe_r^*$, which has
 total volume $1$ on $T$, 
and on $\frakh$ we take
$$\boe_1^*\wedge\ldots\wedge\boe_r^*\wedge\bof_1^*\wedge\ldots\wedge\bof_k^*\wedge
\bog_1^*\wedge\ldots\wedge\bog_k^*.$$   This induces the measure 
$\bof_1^*\wedge\ldots\wedge\bof_k^*\wedge
\bog_1^*\wedge\ldots\wedge\bog_k^*$ on $\frakh/\frakk = \frakh/\ker
\rho$.  

We now use the anchor map $\rho$ (or,
more precisely, ${\rho^{-1}}^*$), to transfer the last measure to the
tangent space $\rho (\frakh)$ to the adjoint orbit through a typical point
$\boe$ (which may be taken in $\frakt$ since the latter intersects
every orbit).
Since, at the basepoint $\boe\in \frakt$, $\rho$ is just
$[\boe,\cdot]$, so that ${\rho^{-1}}^* (\bof_j \wedge \bog_j) = (\sigma_j
(\boe))^{-2} \bof_j \wedge \bog_j,$
we obtain the measure
\[
\prod_j (\sigma_j(\boe))^{-2} \bof_1^*\wedge\ldots\wedge\bof_k^*\wedge
\bog_1^*\wedge\ldots\wedge\bog_k^*
\]
 on the tangent space to the
adjoint orbit.  Dividing the given measure on $\frakh$ by this one, we
obtain 

\begin{equation}
\label{eq-adjointorbit}
\prod_j (\sigma_j(\boe))^{2}
\boe_1^*\wedge\ldots\wedge\boe_r^*
\end{equation}
 on the normal space to the orbit,
i.e. on the tangent space to $\frakh/H$, identified with a tangent
space to $\frakt$.   
 
Since the volume of each isotropy group in the chosen measure is equal
to $1$, the expression (\ref{eq-adjointorbit}) is also the induced
measure on the stack $\frakh//H$.  We may compare this expression with
a calculation in \cite{do-wi:harmonic}.  According to the formula
there on the bottom of p.~22, for any function $\phi$ on $\frakh/H$, 
\begin{equation}
\int_\frakh \phi (p(X))dX = \int_{\frakt^+}
\prod_j(\sigma_j(\boe))^{2}\phi(\boe)d\boe, 
\end{equation}
where $\frakh/H$ is identified with the positive Weyl chamber
$\frakt^+$, $p:\frakh\to \frakh/H$ is the natural projection, $dX$ is the
measure on $\frakh$ which agrees at $0$ (via the exponential map) with the invariant measure on
$H$ with total volume $1$, and $d\boe = \boe_1^*\wedge\ldots\wedge\boe_r^*.$
Thus, our ``natural'' measure on $\frakh//H$, which involved no choices, may be identified with
the push-forward under $p$ of the measure $dX$ on $\frakh$ normalized
as described above.  

If the group $H$ is noncompact, but is 
semisimple and reductive, its adjoint representation still defines 
a proper action when restricted to
the open subset $\cald_H \subset \frakh$ of strongly stable
elements.  These are defined in \cite{we:discrete}, where
it is shown that the adjoint action of $H$ on $\cald_H$ is proper.
Several characterizations of strongly stable elements are given there;
the shortest of these to state is that the
strongly stable elements are those which
belong to the Lie algebra of a unique maximal compact subgroup of
$K$.  It is also shown there (in slightly different terms) that the action
groupoid $H \times \cald_H\arrows\cald_H$ is equivalent to its
restriction $K \times \cale_K\arrows\cale_K$, where $K$ is any maximal
compact subgroup of $H$ and $\cale_K = \cald_H \cap \frakk$.  (The
elements $v$ of $\cale_K$ are characterized by the condition that 
the vector field $\mu_S$ on the symmetric space $S=H/K$ given by
    the infinitesimal action of $\mu \in \calh $ has a nondegenerate zero at the
    coset $eK$.)  
Thus, the natural measure on $\cald_H//H$ is the same as that on
$\cale_K//K$.  

A similar reduction is possible for the conjugation action of $H$ on
itself, following the analysis of Demazure \cite{de:krein}.

\subsection{Symplectic groupoids and Poisson manifolds}
A Lie groupoid $G\arrows G_0$ is a {\bf symplectic groupoid}  if $G$ is equipped
with a symplectic structure $\omega$ which is multiplicative in the
sense that $m^*\omega = \omega_1 + \omega_2$ on the space $G_2\subset
G \times G$ of
composable pairs, where $m:G_2 \to G$ is the product operation, and
$\omega_1$ and $\omega_2$ are the pull-backs of $\omega$ by the first
and second projections from $G_2$ to $G$.  Some standard facts (see, for
example \cite{va:lectures}) about symplectic groupoids are: $G_0$ carries
a unique Poisson structure $\Pi$ for which $l$ and $r$ are
Poisson and anti-Poisson maps respectively; the units form a
lagrangian submanifold of $G$; the Lie algebroid of $G$ is naturally
isomorphic to $T^*G_0$ with the anchor map $T^*G_0\to TG_0$ given by $\Pi$ and
a bracket for which $\{df,dg\}=d\{f,g\}$ for all functions $f$ and $g$
on $G_0$.  Conversely, if $G_0$ is any Poisson manifold for which the
associated Lie algebroid structure on $T^*G_0$ is integrable to a
groupoid, then the integrating groupoid with simply-connected source
fibres is a symplectic groupoid with underlying Poisson manifold
$G_0$.  If we allow integration by stacks rather than just manifolds,
then there is a bijective (up to natural isomorphisms) correspondence
between Poisson manifolds and symplectic groupoids with
simply-connected source fibres \cite{ts-zh:poisson}.  

Writing $A$ as usual for the Lie algebroid $T^*G_0$, we have a natural
isomorphism between $Q_A$ and the tensor square $(\wedgetop
T^*G_0)^{\otimes 2}.$   $Q_A$ therefore carries a natural orientation,
and its positive sections are the squares of nowhere-vanishing
sections $\nu$ of $\wedgetop T^*G_0$.  Assuming $G$ to have connected
source fibres, we have (whether or not $\nu$ is nowhere-vanishing) 
that $\nu^2$ is $G$-invariant exactly when $\nu$ is invariant under
all hamiltonian vector fields.  

\subsection{Symplectic manifolds}

Let $(S,\omega_S)$ be a $2m$-dimensional, 
 connected\footnote{From the point of view of Poisson geometry, we should
  perhaps call a manifold with nondegenerate closed 2-form
 ``symplectic'' {\em only} if it is connected, since
  otherwise it has more than one symplectic leaf!} 
 symplectic manifold, considered as a Poisson manifold.
The source-connected symplectic groupoids for $S$ (all of them
transitive) are just the
quotients of the fundamental groupoid $\pi(S)$ associated with normal
subgroups of the fundamental group of $S$.  For any such groupoid, 
the invariant sections of $Q_A$ are just the constant multiples $\lambda=c\nu_S^2$
of the square of the Liouville measure $\nu_S = \frac{(-1)^{m(m-1)/2}}{m!}\omega_S^m.$  

Now let $G\arrows G_0$ be the symplectic groupoid of $S$ associated
with the subgroup $K$ of the fundamental group.  This 
groupoid is proper
when $K$ is finite; two natural choices are the trivial group 
and the fundamental group itself, when it is finite.  It is equivalent
to the groupoid $K\arrows {\rm pt}$; i.e. $G_0/G$ is just a point,
and $G_0//G = BK$.  We recall from the Introduction that
$\#(BK)=1/\#(K)$ (in particular, we get the value $1$ when $K$ is the
trivial group and $G$ is the pair groupoid), but we are interested
here in the volume $\vol (BK)$, which will depend on the choice of 
the constant $c$ in $\lambda$.  To compute this volume, we begin by
factoring $c\nu_s^2$ as the product of the sections $ a^{-1}=c\nu_s$ of
$\wedgetop A^*$ and $b = \nu_S$ of $\wedgetop T^* G_0$.  Next, we
choose the section $1$ of the trivial bundle $\ker \rho$.  Noting
first that integrating $1$ over the isotropy groups will give $\#(K)$, 
we next observe that the induced section of $\wedgetop(A/\ker\rho)^*$ is again
$a^{-1} = c\nu_S.$  Its inverse is the section $a=c^{-1} \nu_S^{-1}$ of
$\wedgetop (A/\ker\rho)=\wedgetop (T^*G_0)^*.$  
Now we must transfer this to a section of $\wedgetop
T^*G_0$ by using ${\rho^{-1}}^*$,  
with $\rho:T^*S\to TS$ given by the Poisson structure inverse to
$\omega_S$.  It is not hard to see (for instance by using 
symplectic bases) that $\rho^*$ pulls back $\nu_S$ to $\nu_S^{-1}$, so 
${\rho^{-1}}^*$ transfers $c^{-1}\nu_S^{-1}$ to $c^{-1}\nu_S$.
Finally, we divide this measure along the ``orbits'' into the measure
$\nu_S$ on the base of the groupoid to obtain the measure $c$ on the 
orbit space.  Since this orbit space consists of a single point, we
conclude that $\vol_{c \nu_S^2}(G_0//G) = c/\#(K)$.  Observing that this agrees
with the cardinality just when $c = 1$, we are led to the
conclusion that the Liouville measure (including the factor of $1/n!$)
is really the natural one on a symplectic manifold.   (Of course, this
is just a consequence of the fact that the Liouville measure is the
only one, up to sign, which corresponds to its inverse by the Poisson
structure.)  

\subsection{Regular Poisson manifolds}
We look here at a very simple example.
  Let $G_0=P=\Sigma\times \reals$ be the
Poisson manifold given by a family $\omega_t$ of symplectic
structures on a 2-sphere $\Sigma$, parametrized by $t\in
\reals$.  If the area $V(t) = \int_\Sigma \omega_t$ of the symplectic
leaf $\Sigma \times \{t\}$ has no critical points
as a function of $t$, then this Poisson manifold is integrable.  
(See, for example, Sections 5 and 7 of \cite{da:groupoides}).  Its
symplectic groupoid $G$ is a circle bundle over $\Sigma \times \Sigma
\times \reals$ for which the map $(l,r)$ takes the entire fibre over
$(x,y,t)$ to $((x,t),(y,t))$.  The fibres over the diagonal points
$(x,t),(x,t)$ are the isotropy groups, and the corresponding Lie
algebra bundle is naturally isomorphic to $T^*\reals$.  (For any
Poisson manifold at a regular point, the isotropy algebras of a
symplectic groupoid are naturally isomorphic to the conormal spaces of
the symplectic leaves.)  Under this isomorphism, the kernel of the
exponential map from the isotropy algebra at $(x,t)$ to the
corresponding isotropy group consists of the integer multiples of $dV
= V'(t)dt$. 

Let us now choose the section $\lambda = f(t)(\omega_t \wedge dt)^2$
of $Q_A$ and compute the corresponding measure on the stack $G_0//G$,
which is a bundle of $BS^1$'s over $\reals$.  As we did in the
symplectic case, we factor $\lambda$ as the product of the sections 
$a^{-1} = f(t)\omega_t\wedge dt$ of $\wedgetop A^*$ and
$b=\omega_t\wedge dt$ of $\wedgetop T^*G_0$.  Writing $\tau$ for the
coordinate on the isotropy groups corresponding to $dt$, we choose the
measure $d\tau$, for which the volume of the isotropy group at
$(x,t)$ is equal to $V'(t)$.  
The induced measure on $A/\ker \rho$ is then $f(t)^{-1}(\omega_t)^{-1}$, which
pushes forward under ${\rho^{-1}}^*$ to the measure
$f(t)^{-1}\omega_t$ 
along the symplectic leaves.  Dividing this into
the measure $\omega_t \wedge dt$ on the base of our groupoid gives the
measure $f(t)dt$ on the quotient space.  Using now the fact that the
volume of the isotropy group at $t$ is $V'(t)$, we find that the
measure on the stack $G_0//G$ is equal to $(f(t)/V'(t))dt$.  To
better understand this result, we assume $f$ positive, so
that $\lambda$ is the square of the measure $(f(t))^{1/2}\omega_t\wedge dt$ on
$P$.  Dividing the latter by the Liouville measure along the leaves gives
$(f(t))^{1/2} dt$ as ``quotient'' measure on $\reals$.  It agrees
with the stack measure just when $\lambda$ is the square of
$\omega_t\wedge V'(t)dt =\omega_t \wedge dV(t).$  

We may interpret the last calculations as meaning that, just as
Liouville measure is a natural measure on a symplectic manifold,
so $dV$ is a ``natural'' measure on the leaf space of
the Poisson manifold above.  A possible generalization to
arbitrary regular Poisson manifolds goes as follows.  

Let $G_0=P$ be a Poisson manifold which is strongly regular in the
sense that its symplectic leaves are the fibres of a smooth fibration
$P\to M$.  Following \cite{da:groupoides}, we introduce the lattice
(``r\'eseau'') $\Lambda \subset T^*M$ of differentials of periods of the
symplectic forms on the leaves.  It is a (not necessarily closed)
lagrangian submanifold of $T^*M$ for which the projection to $M$ is a
local diffeomorphism and whose intersection with each fibre is a (discrete)
additive subgroup.  The Poisson manifold is integrable to a symplectic
groupoid when $\Lambda$ is closed and a covering space of $M$.  

We will assume further, for simplicity, that the leaves are simply
connected and that the intersection of $\Lambda$ with each cotangent
space is a full lattice, in which case the quotient $T^*M/\Lambda$ is a
bundle of tori whose pull-back to $P$ is the isotropy subgroupoid of
the symplectic groupoid $G$ of $P$.   (The assumption that the leaves are
simply connected implies that the symplectic groupoid is unique up to
isomorphism.)  In this case, the lattice defines an integrable
$GL(n,\mathbb Z)$ structure on $M$ which determines (up to sign) a
natural measure.
 
\begin{conj}  
With notation as above, if the section $\lambda$ of $Q_A$ is the square
of an invariant measure on $P$, factored as the product of the
Liouville measure along the symplectic leaves and (the pull-back of) a
measure $\beta$ on the leaf space $M$, then the
induced measure on the stack $P//G$ agrees with $\beta$ if and only
$\beta$ is the measure associated to the integer affine structure on
$M$.
\end{conj}

\subsection{Duals of Lie algebras}  
If $H$ is a compact Lie group, its coadjoint representation on $\frakh^*$ is
equivalent, via a bi-invariant metric, to its adjoint representation
on $\frakh$.  Different choices of the metric will lead to different
identifications, hence the canonical section of $Q_A$ for the adjoint
action does not lead to a canonical section of $Q_A$ for the coadjoint
action. 

In fact, the natural structure on $\frakh^*$ is its
Lie-Poisson structure, for which $T^*G$ is a symplectic groupoid.  
On an open dense subset, this Poisson
structure is regular and satisfies the hypotheses of the conjecture at
the end of the previous section.  For the special case of $SU(2)$, the
Poisson structure on the complement of the origin has concentric
spheres as its symplectic leaves, and we may choose a linear radial
coordinate $t$ so that the symplectic area of the sphere of radius $t$
is $V(t)=4\pi t$.  The ``natural'' measure on the leaf space is then
$4 \pi dt$, whose product with the Liouville measure along the leaves
is $\omega_t \wedge 4\pi dt = 4\pi t \omega_1 \wedge dt.$  Notice
that this is not a translation-invariant measure!

\end{document}